\documentclass[graybox]{svmult}

\usepackage{mathptmx}       
\usepackage{helvet}         
\usepackage{courier}        
\usepackage{type1cm}        
\usepackage{makeidx}         
\usepackage{graphicx}        
\usepackage{multicol}        
\usepackage[bottom]{footmisc}

\usepackage{amsmath}
\usepackage{amssymb}

\newcommand{\proves}{\vdash}

\makeindex             

\begin{document}

\title*{A New Foundational Crisis in Mathematics, Is It Really Happening?}
\author{Mirna D\v zamonja}
\institute{Mirna D\v zamonja \at School of Mathematics, University of East Anglia, Norwich NR47TJ, UK,  \email{m.dzamonja@uea.ac.uk}
\and equally associated to \at IHPST, Universit\'e Panth\'eon-Sorbonne, Paris, France }
%
%

\maketitle

\abstract{The article reconsiders the position of the foundations of mathematics after the 
discovery of the homotopy type theory HoTT. Discussion that this discovery has generated in the community of mathematicians, philosophers and computer scientists might indicate a new crisis in the foundation of mathematics. By examining the mathematical facts behind HoTT and their relation with the existing foundations, we conclude that the present crisis is not one. We reiterate a pluralist vision of the foundations of mathematics.\\
The article contains a short survey of the mathematical and historical background needed to understand the main tenets of the foundational issues.}

\section{Introduction}
\label{sec:intro} It seems by all evidence about the differences between set theory and homotopy type theory HoTT that there is a crisis in the foundations of mathematics. However, neither mathematicians nor even the philosophers have spoken enough about it. The preceeding crisis, happening at the beginning of the 20th century and sparked by the Russell paradox, created many animated discussions and engaged a large number of mathematicians- not necessarily `logicians' (as mathematical logic as a separate field of mathematics had not yet existed, properly speaking).
All the mathematicians were interested in foundations and so were many philosophers.  That was the time when mathematics formed a tameable, unified, body of research, and mathematicians could easily understand each other. The folklore has it that Hilbert was the last mathematician who knew all of mathematics, 
but he was not the only one: let us mention Poincar\'e, Picard and Von Neumann, for example. The discussion about foundation of mathematics engaged a wide audience of mathematicians.
Today, logicians tend to hold foundational discussions with very little success in engaging the mathematical community. For example, mathematicians spent 358 years solving Fermat's last problem. The solution, due to Wiles \cite{Wiles} could not have possibly been the one that Fermat meant to write in his famous margin in 1637, as it involves mathematical techniques that did not exist even in a trace form at the time of Fermat. Could Fermat have found a different proof? If we prove that Fermat's last problem does not have a solution in Peano Arithmetic PA, which encompasses all the tools available to Fermat, then it is pretty clear that Fermat could not have solved it. Some FOM (Foundations of Mathematics) discussion groups seem convinced that this is the case, 
yet Angus Macintyre seems to have proved that Wiles' proof actually can be carried out in PA
\cite{Angus}. However the general mathematical audience has not yet noticed this `rumination by the logicians', as much as everybody would really love to know if Fermat had actually proved 
it. 

The truth is that mathematical logic has become a very specialised subject, considered by mathematicians as an esoteric subject of little interest to their daily life. Therefore, even the mathematical foundations, since they are immersed in logic, are considered as esoteric. (It must be admitted that this feeling in the mathematical community is not completely decoupled from the rather dry writing style that logic papers tend to have). The purpose of this article is, firstly, to recall that just exactly the opposite is true, thinking about foundations is not esoteric and, if there is indeed a crisis in the foundation of mathematics, then {\em everybody} should be concerned. Not only the logicians and the mathematicians but also all the scientists, who base their reasoning upon logic and mathematics. One cannot build a high-rise without knowing how stable are the pillars.\footnote{This metaphor might be overused. But I still like it.}  Secondly, in the main body of the article we shall describe the foundational question under consideration in such terms that every mathematician, philosopher of mathematics or computer scientist can follow the exposition, with the purpose of engaging a wider community in the foundational discussion. Lastly, we shall give our own opinion, calming the fire of panic and arguing that the crisis is relative and solvable and that it brings us more of the positive than of the negative information. We shall argue for a vision of foundations that envelops this information into a new stage of the development of the existing foundations of mathematics.

\section{For a Mathematician}\label{sec:maths} Many of our glorious ancestors in set theory in the 1960s and later put an enormous effort to show to mathematicians that the fact that there were simple statements, such as the continuum hypothesis\footnote{stating that every infinite subset of the set of the real numbers is bijective either with the set of the natural numbers or with the set of the real numbers}, about simple objects, such as the set of the real numbers, that turned out to be independent of the axioms of set theory - then almost universally considered as {\em the}
foundation of mathematics - should be of interest to every mathematician on the planet.  This effort does not seem to have been very successful, even though the independence can appear almost everywhere, see \cite{Fremlin} for some early examples in various areas of mathematics. Well, the point is that it is {\em  almost} everywhere, and many mathematicians  have had the reaction to withdraw to those areas of mathematics that were immune to independence and to relegate any independence to the inferior rang of `pathological' objects. Infinite combinatorics- pathological, uncountable measures- pathological, non-separable Banach spaces- pathological, independence- ``can't possibly happen to me". 

Fine, but what if somebody tells you that your art of mathematics is going to go the way of the games chess and go, it is going to be replaced by formal arguments doable by a machine which can do it better than you or any other human? That is how the specific program of formalisation of mathematics through univalent foundations teamed with automatic theorem proving, presented by a winning mathematician who invented it, looks to a naked eye. 

Ironically, it is exactly those areas of mathematics most remote from independence, that are protected from independence by their inherent definability and constructibility, which are in danger of this new formalisation and computerisation.  Formalisation of mathematics most often uses constructive mathematics. `Pathological' objects in general are not constructive, they often rely on the axiom of choice and further set-theoretic machinery, with no specific witness for existential properties. This is what makes them prone to independence, but rather safe from formalisation\footnote{although see later on Isabelle and see Mizar's project website http://mizar.org/project/}.  So,  one way or another, foundations have not been kind to mathematicians, and it looks like one should start paying attention and ask the question if the mathematics we are producing today is safe from a foundational and existential crisis. 

Ever since the 19th century most (or just many?) people have considered set theory as the foundations of mathematics. This is in spite of the great crisis in the foundations of mathematics brought about by Russell's paradox in 1902 and the subsequent introduction to mathematics not only of several different axiomatisations of set theory\footnote{In this paper we take Zermelo-Fraenkel set theory with choice (ZFC) for our basis and we freely dispose of a sufficient amount of large cardinals.} and various very successful alternative approaches to foundations of mathematics, such as category theory, type theory and the theory of toposes. While it is clear, for example, that it is much more reasonable to do certain parts of mathematics by using category theory, it has been known that category theory and set theory with large cardinals are to some extent bi-interpretable (see e.g. \cite{Feferman}) and apart from philosophical reasons, there is no reason for an `ordinary mathematician' to take his time from his mathematical considerations in order to think if one or another foundation fits better with the purpose at hand.

In 2006 a new kid appeared on the block. It is the homotopy type theory including univalent foundations.
These foundations are different from set-theoretic foundations. In some respect these foundations appear to contradict the set-theoretic foundations. A great disquiet, or perhaps excitement about the newly proposed foundations can be seen by the following quote by
Michael Harris\footnote{a mathematician who does care about the foundations of mathematics} in \cite{Harris}:
\begin{quotation} 
It's impossible to overstate the consequences for philosophy, especially the philosophy of Mathematics, if Voevodsky's proposed new Foundations were adopted.
\end{quotation} 

Our thesis is that there is no contradiction in sight, moreover, these foundations complement the set-theoretic foundations. We feel that any suggestion of a contradiction between the two or the idea that either should or could be replaced by the other comes more from a lack of understanding than anything else. We attempt in this article to give a simple account of the most relevant developments
in order to demonstrate our thesis. We shall lead the reader through a presentation of the basics of several different parts of mathematics and logic that have gone into the developments in question and sometimes we shall have to simplify to a great extent in order to make ourselves understandable. The reader interested in the details will easily be able to find them, either in the other chapters of this book or elsewhere in the literature and references mentioned here. 

\section{On Set Theory and Type Theories}\label{setsandtypes} Set theory is based on the classical first order logic. Basic entities are sets and they are completely determined by their elements: the Axiom of 
Extensionality states that for any sets $A$ and $B$,
\[
[\forall x (x\in A\iff x\in B)]\iff A=B\mbox{.     }\footnote{The intended meaning of the symbol $\in$ is that it stands for the actual real' membership relation. A common misunderstanding is that this is the case with all models of set theory. However, the setting of set theory in the logic of first order, does not allow to rule out models in which the symbol is interpreted by other means. They have to be dealt by separate means, including the Mostowski collapse theorem.}
\]
The Law of Excluded Middle holds and the Axiom of Choice is assumed (as we are working with ZFC). This type of set theory was formalised by Ernst Zermelo \cite{Zermelo} and ramified in the work of Zermelo and Fraenkel, as well as Skolem. A detailed history of the development of ZFC and some other set theories can be found in \cite{FHL}. 

Univalent foundations are a vision of Vladimir Voevodsky of a 
foundational system for mathematics in which the basic objects are homotopy types. They obey a type theory satisfying the univalence axiom. They are formalisable in a computer proof assistant and are constructivist. To introduce univalent foundations, let us start by discussing type theory. We present the types together with the idea of a proof, which is not how these concepts were developed initially. 

In set theory, basic variables are interpreted as sets and atomic formulas are of the form $x=y$ and $x\in y$. Naive set theory in which every collection of the form $\{x:\,\varphi(x)\}$ forms a set, leads to paradoxes, the most celebrated of which is the one developed by Russell in 1902\footnote{Zermelo seems to have discovered the paradox about at the same time, although he did not publish it.}. Presuming that $B=\{x: x\notin x\}$ is a set, we see that the 
question if $B\in B$ does not have a unique answer. Set theory resolved this issue by various axiomatisations, in our case of ZFC by introduction of Restricted Comprehension by Zermelo in
1908. This type of Comprehension only allows formation of sets of the form $\{x \in y :\,\varphi(x)\}$, that is, the comprehension of all elements of a given set with a given property, into a set. Russell's own solution to the paradox was to introduce type theory. There are many versions of 
type theory, but the following version from \cite{FHL} is equivalent to the one used by Russell and Whitehead in \cite{Whitehead.1990} in which they succeeded to give foundations of the entire extent of the mathematics then known in terms of type theory.

In this context, each variable has a type, which is a natural number, and writing $x^i$ means that $x$ is a variable of type $i$. Atomic formulas are of the form $x^i \in y^{i+1}$ and $x^i =
z^i$. That is, a variable can only be equal to a variable of the same type and can only belong to a variable of type exactly one higher. The Axiom of Infinity is added to level 0 and Axiom of Choice at every level. 

The advantageous formalisation power of type theory is to some extent balanced off by the fact that is rather cumbersome to write mathematics in this way. The reader may find it entertaining to contemplate the notions of complement or cardinality in this setting.

For the homotopy type theory it is more useful to think of the Dependent Type Theory, developed by Martin-L\"of starting from the 1960s.\footnote{In our explanation, we shall privilege intuition to exactness. A much more precise, yet still rather compact, survey of Type Theory and Univalent Foundations can be found in Thierry Coquand's article for the Stanford Encyclopedia of Philosophy \cite{CoquandStanford}. A longer survey of these matters is in the recent article by
Daniel Grayson \cite{Grayson}. Another recent article, by Martin Escardo \cite{Escardo} gives an entirely self-sufficient presentation of the Axiom of Univalence. We have partially followed an article by Michael Shulman and its notation,
\cite{Schulman}.}
Although the original idea of dependent types evolves from Martin-L\"of's study of randomness in probability,  we can see dependent type theory in a computational context, like a computer language (this idea goes back to de Bruijn). It is particularly well suited for the theory of proofs, as we shall now see. 

Type theory does not have semantics, only syntax.
The basic expressions are of the form {\bf term:Type}, for example $x:{\mathbb N}$ and
$y:{\mathbb R}$, which mean that $x$ is a variable of the non-negative integer type and $y$ a variable of the
real type. Such expressions, variable typings, are called {\em contexts}.
From this we can advance to
{\em proofs} or {\em judgements} using some predetermined rules. For example
\[
x:{\mathbb N}, y:{\mathbb R}\proves x \times y: {\mathbb N}\times {\mathbb R}.
\]
The symbol $\proves$ means {`}proves'. In this judgement we can see a typical operation on types, which in this case is the product of two. A list of all basic types and the operations on them 
can be found on the slides by Coquand \cite{CoquandBourbaki}.
It is useful to think of this as a category, in which objects are the contexts and the morphisms are the judgements. This is called the {\em Classifying Category} and denoted {\bf Ctx}. So, when we write 
\[
\Gamma\proves a:A
\]
we consider $a$ as a morphism from $\Gamma$ to $A$. Contexts with variables are `types' and those without variables, where we substituted a variable by a constant, are morphisms.
It is exactly the way to handle the idea of substitution which makes this into the {\em dependent type theory}. For example, we allow families of contexts, such as $B$ {\bf type}, where we interpret
\[
\Gamma\proves B \mbox{ {\bf  type}}
\]
as a judgement allowing to conclude from $\Gamma$ that there is an object of type $B$ (or $B$ is a type in $\Gamma$). The
intuition is that $B$ `depends' on $\Gamma$. The dependent types also can be combined, for example
\[
{n:\mathbb N}\proves C_n \mbox{ {\bf  type}}
\]
corresponds to ``for all $n$ we have an object of type $C_n$''.\footnote{Of course, type theory has no internal logic in the sense of semantics, but the operations of product and join correspond to the existential and universal quantifiers, under the Curry-Howard correspondence, described below. See \cite{Escardo}. The intuitionistic higher-order logic interpretable in MLTT is often referred to as the ‘internal logic of (dependent) type theory.}. We can think of the dependencies as factor or slice categories, where for example the statement $\Gamma\proves B$ {\bf type} corresponds to the existence of type $B$ in the category {\bf Ctx}$/\Gamma$, the category of
all contexts in which $\Gamma$ holds, or precisely to the statement that $B$ is a type in that category. Shulman \cite{Schulman} states: ``The point is that substitution into a dependent type presents the pullback functor between slice categories". Intuitively speaking this means that substitutions behaves as we are used to having it,
where 'free' variables can be replaced by other free variables of the same type.

These developments lead to a proof system by allowing deductions, which are the axioms of judgements. These axioms come in two kinds: ones that allow us to construct a type 
and others which allow to eliminate in order to prove equivalence. This resembles, for example, the natural deduction, but in this context it is more complex. An example of a deduction rule looks like this:  
\[
\frac{\Gamma\proves A \mbox{ {\bf  type}}; \Gamma; n:A\proves B \mbox{ {\bf  type}}}
{\Gamma\proves \prod_{n:A} B \mbox{ {\bf  type}}}
\]
where the premise is written in the numerator and the conclusion in the denominator. 
This model of reasoning is very formal and well suited for implementation on a computer, hence proof assistants such as Coq. Choosing a particular collection of type constructors specifies the rules of a particular type theory, that is, dependent type theory is not a unique theory but rather a template for producing proof theories. 

It is possible to interpret propositional logic in type theory, namely we interpret propositions $P$ as types with at most one element. $P$ is `true' if there is an element in it (these are types of level 1). From the point of view of proofs, saying that type is inhabited amounts to saying that there is a proof of that type. We can also interpret sets, as properties (types of level 2)\footnote{This may look like type theory is a generalisation of set theory, as mentioned in some references, but this is not true since set theory is not only made of objects but also of the underlying logic and the axioms.}. By a property of elements of a type $A$ we mean a judgement of the form
\[
                         x:A\proves P:\Omega.
\]     
Under this interpretation, we obtain the Curry-Howard correspondence \cite{Howard} between the operations of intuitionistic propositional logic and (some of) the constructors of type theory. So theorems of intuitionistic propositional logic are provable in the deduction system.  The Law of Excluded Middle LEM does not fall into this correspondence. LEM is addable artificially as a judgement, but then it tends to cause difficulties in the execution of the proof assistants, although it seems relatively harmless in Coq.\footnote{Andrej Bauer in private correspondence (27 February 2018) says ``Coq does not loop if we postulate excluded middle, it just gets stuck on some computations, and then requires human intervention to get unstuck."}
There are proof assistants based on higher order logic, such as Isabelle/HOL and Isabelle/ZF which are able to deal with AC and LEM. For example, Lawrence C. Paulson in \cite{Paulson} mechanised G{\"o}del's proof that the consistency of ZF implies that of ZFC using the proof assistant Isabelle. The HoTT libraries of computer proofs are based on Coq.

We note that the full Axiom of Choice is inamicable to the intuitionistic logic in the sense that by proofs of Diaconescu \cite{Diaconescu} and Goodman and Myhill \cite{GoodmanMyhill},  the Axiom of Choice implies the Law of Excluded Middle. In a model of  dependent type theory, the Axiom of Choice can be enforced on the level of sets, at the likely price of losing the connection with the proof assistants. A very interesting article about the relation of type theory and the Axiom of Choice is Martin-L{\"o}f's \cite{100Zermelo}, in which it is shown that when a certain amount of choice is added to type theory the resulting theory interprets ZFC.

\section{An Identity Crisis and Univalence}\label{identity} To explain what the univalence is, let us first discuss some intuitions.
Frege in the celebrated Appendix of his \cite{Frege} (in reaction to Russell's paradox and attempts to solve it) stated that
\begin{quotation}
Identity is a relation given to us in such a specific form that it is inconceivable that various kinds of it should occur.
\end{quotation} 
Indeed, in set theory, there is only one kind of equality, given by the Axiom of Extensionality and in model theory it is well known that we can always assume that the equality
relation in any model is interpreted by the true equality. 

Martin L{\" o}f's  type theory MLTT builds up on Dependent Type Theory by introducing some new types, including universes $U_0, U_1, \ldots$ such that $U_n: U_{n+1}$ and the identity type $Id_A(x,y)$ for propositions $x,y$. This type is the type of proofs of equality between $x$ and $y$, which are both of type $A$. So if this type is inhabited, it means that the elements $x$ and $y$ of type $A$ are provably equivalent. Another type, called ${\rm refl}$ proves that any $x$ is equal to itself. This setting gives rise to {\em two} different kinds of identity: definitional $A=B$ (for types) and $x = y : A$ for elements of a type, and the propositional one $Id_A(x,y)$. Furthermore, building on $Id_A(x,y)$ we can also declare that any two proofs for this identity are equivalent etc., extending the propositional equivalence to all levels of types. Then it is clear (in fact it follows from the rule of identity introduction)  that the definitional identity implies the propositional one, but it is not at all clear that the converse should be true. 

This problem, which is to some sense the completeness problem for proof theory, was noticed early on by Martin-L\"of and others, who tried to resolve it by inserting explicit judgements to this extent. However, in the full generality,  this leads to paradoxes.
The sameness of the two notions of equality became a known problem and search was open for more subtle models of dependent type theory where this sameness is true.\footnote{This leads to intensional versus extensional type theory and is a long and interesting story, see \cite{Hofmannthesis} for more.}

Intuitively speaking, Voevodsky's univalence fullfils this dream. It says that the identity type 
$\mbox{Id}(X, Y )$ is precisely the type $\mbox{Eq}(X, Y )$ of equivalences (which we have not defined), in the sense of being in one- to-one correspondence with it.
To quote \cite{Grayson}, 
{\em it offers a language for mathematics invariant under ` equivalence'}. This intuitive understanding of the axiom, present in most presentations of univalence for general audience (hence this one), is somewhat misleading. Much is to be gained by understanding a more exact definition of the univalence, which is possible to do by following the short and exactly-to-the point \cite{Escardo}.  
The conclusion is that univalence is a property of Martin-L{\" o}f's identity type of a universe of types. Without univalence, ${\rm refl}$ is the only way to construct elements of the identity type, with univalence this is also possible by other means.
The {\em univalence axiom} states that a certain type is inhabited, the univalence type. It takes a number of steps to build the univalence type, for these steps the reader may consult \S3 of
\cite{Escardo}.

Voevodsky's HoTT is an extension of MLTT by the Univalence Axiom, using topological insights which we shall describe in \S5. It should be noted that  the Univalence Axiom is not a consequence of MLTT as there is also a model of MLTT where it fails, namely any model of Streicher's Axiom K (see \cite{Streicher}). 

\section{The Topological View}\label{topology} Vladimir Voevodsky came into the field in a rather round-about way, and one may say, somewhat naively. Paraphrasing his talk in the 2013 Logic Colloquium in Evora, Portugal, we learn that after obtaining a Fields medal in 2002, he was feeling increasingly isolated in the position to publish any proof he wanted (due to his reputation) and to have a shrinking audience of real readers (due to the complexity of the proofs in question). He realised that he was not alone in the situation and judged that mathematics badly needed - computerised proof verification. Only then did he understand that there already existed a vibrant community of researchers working not only on proof verification but also on proof assistants. At learning of the relevant issues, Voevodsky realised that the univalency question had appeared in a different part of mathematics under a different name and that it had been successfully solved. The field was homotopy theory and the key to the solution is the notion of an $\infty$-groupoid (a notion that was already used by some other researchers in the logic community, as it happens, see below for historical remarks). The main idea of Voevodsky was to interpret the category {\bf Ctx} as the homotopy category, so that types become spaces and morphisms. Propositional equality then becomes the homotopic equivalence and the structure obtained becomes an $\infty$-groupoid. Finally, Voevodsky's vision become the foundation
of mathematics expressed in the language of type theory in which the objects are represented in
an $\infty$-groupoid. This vision is succinctly represented in Coquand's Bourbaki seminar presentation \cite{CoquandBourbaki}. 

For us here, let us briefly recall the relevant topological notions and their history, discovered and developed in 1950s by Daniel Kan, building on an even earlier work by Henry Whitehead.
We start with the {\em simplex category} $\Delta$, whose objects are linear orders on a
natural number $n=\{0, 1, 2, \ldots, n-1\}$ and whose morphisms are (non-strictly) order-preserving functions.   We fix some large enough infinite cardinal number (may be a large cardinal) and consider all sets that have cardinality smaller than that cardinal. This category is called {\bf Set}. A {\em simplicial set} $X$ is a contravariant functor $X:\,\Delta\to$
{\bf Set}. Simplicial sets form a category under natural (the word natural here has a technical meaning) transformations and this category is called {\bf sSet}.

The point is that it is a classical result in algebraic topology that a certain fibration {\bf sSet}$/W$ of {\bf sSet} called Kan complexes is equivalent to the category of `nice' topological spaces and homotopy classes of maps. Because of this fact one can, when working up to homotopy, think of simplicial sets as of combinatorial representations of shapes or spaces, of simplicial paths as paths on these spaces etc. 
Kan complexes form a model of an $\infty$-groupoid, the formal definition of which requires higher category theory, but which in essence deals with a category with objects of every type $n$ for $n<\omega$ and $n$-morphisms. The homotopy types of the geometric realizations of Kan complexes give models for every homotopy type of spaces with only one non-trivial homotopy group, the so called Eilenberg-Mac Lane spaces. It is conjectured that there are many different equivalent models for $\infty$-groupoids all which can be realized as homotopy types, some of which are known. 

Grothendick's {\em homotopy hypothesis} states that the $\infty$-groupoids are equivalent 
as in the equivalence of categories, to 
topological spaces.

A central theorem of the homotopy type theory with univalence, and one which gave it its name is the following:

\begin{theorem}[Voevodsky \cite{KapulkinLumsdaine}] \label{Voy} Modulo the existence of two inaccessible cardinals, it is consistent that {\bf sSet}$/W$ forms a model of Martin L{\" o}f's type theory with the Univalence Axiom.
\end{theorem}

In the proof one can also assume LEM on the level of propositions and the Axiom of Choice on the level of sets. Working in this particular model is what became known as the {\em univalent foundations}, although some other models were found later, as we shall mention. 

The revolutionary result of Voevodsky had some precursors and parallels, notably in the work of Michael Makkai 
\cite{Makkai}, Hofmann and
Streicher \cite{Streicher} and Steve Awodey (see \cite{Awodey1}); 
Martin Hofmann's Ph.D. thesis \cite{Hofmannthesis} already mentioned above, for example, contains a study of the groupoid of intensional equivalences. 
See \cite{wikihott} for a detailed historical description.
It also has descendants, which are well worth mentioning.  Another model of Univalence Axiom and with less consistency strength, the cubical model, was found by Bezem, Coquand and Huber \cite{Bezem}. Work of Peter Aczel \cite{Aczel} shows that dependent type theory for constructive sets has a much weaker proof theoretic strength than that of ZF. Namely, this model can be done in constructive set theory CZF with universes. This set theory is proof theoretically strictly weaker than ZFC, and has the same strength as MLTT/dependent type theory (without the Univalence Axiom). A corollary of the work by Bezem, Coquand and Huber \cite{Bezem} is that adding the Univalence Axiom does not add any proof theoretic power to dependent type theory. The story of the proof-theoretic strength of these systems is rather involved. A lot of results and the history can be found in Michael Rathjen's  paper
\cite{Rathjen} and, specifically on the strength of the univalence axiom, in the upcoming paper
\cite{Rathjenuni}.

Furthermore, the relation between this and the classical logic and set theory is also very intricate. Work of
Crosilla and Rathjen \cite{CrosillaRathjen} shows that the proof-theoretic strength of CZF+ LEM+universes is at least ZFC+ unboundedly many inaccessibles. So it is LEM that adds strength to univalent foundations. 

 \section{Conclusion, Crisis or No Crisis}\label{sec:3} Much discussion has been
 generated by the advent of the univalent foundations and this volume is certainly a forum for it,
 as is the HoTT book \cite{HoTT}. To certain extreme fans of HoTT, some of them present at the conference ``Foundations of Mathematics: Univalent Foundations and Set Theory'' in  Bielefeld in 2016 which initiated this volume, this  subject means the end of classical mathematics as we know it, as we can just move to a model of univalent foundations and let the ever-developing libraries of computer proofs do all our mathematics for us. A mathematician of the future becomes a new incarnation of a chess player, not capable to come up with anything weakly comparable to what a computer can do on the same subject. The less ardent fans have called for a measured reflexion on the choice between this foundation and the classical ones, implying perhaps the eventual adoption of one or the other. The HoTT book states in the Introduction: ``This suggests a new conception of foundations of mathematics, with intrinsic homotopical content, an “invariant” conception of the objects of mathematics — and convenient machine implementations, which can serve as a practical aid to the working mathematician. This is the Univalent Foundations program. ... we therefore believe that univalent foundations will eventually become a viable alternative to set theory as the “implicit foundation” for the unformalized mathematics done by most mathematicians." Viewed from this angle, the situation looks like quite a serious crisis in foundations, putting two foundational programmes in a competition. 
 
Our view is that this crisis is not one and that, as ever in mathematics, mathematical results speak for themselves. First, there is no obvious reason of why one should choose to limit oneself to a specific model of mathematics, be it the simplicial set model or the cubical one or, for that matter, the constructible universe $L$ or some other set-theoretic universe; the prevalence of consistency results in mathematics has rendered us very careful on this issue. But there are also more tangible observations following from the above results. If a theory has the consistency strength less or equal ZFC+ some amount between 1 and 2 of inaccessible cardinals, then this theory cannot resolve any of the questions for which we know that a larger consistency strength is required. This includes many questions in modern logic, such as projective determinacy, the existence of a measurable cardinal, the failure of the singular cardinal hypothesis and so forth. If we take the view that large cardinals represent the ultimate concept of the `actual infinity' or the 'higher infinite', then these notions cannot be formalised by univalent foundations\footnote{By the `higher infinite' we do not mean the uncountable, which can be formalised in type theory. Rather, we follow the usage the reference book on large cardinals by Kanamori \cite{Kanamori}.}. On another topic, if a theory has a much weaker strength than that same theory when accompanied by the Law of Excluded Middle, then that theory cannot answer all questions which need this law and cannot reproduce all proofs
where this law is used. Some proofs by contradiction have a translation which does not require this law, some can be formalised, but the results quoted above show that this cannot be the case with all such proofs. And much of modern mathematics uses proofs by contradiction. So, there is no way of formalising into a computerised proof in HoTT all that is known or will be known to mathematicians, even those who are not interested in the actual infinity but simply in the use of classical logic. 

Of course, there are and have been many mathematicians (Kronecker, Poincar\'e to some extent and Brouwer to some extent, for example) who have rejected either the actual infinity or the Law of Excluded Middle or both. The fact is that much of what they, in that vison, considered as mathematics is now within the reach of computer formalisation. Is this a bad or a good fact? 
It is neither bad or good. Mathematics is not about judgement but about furthering human knowledge and it recognises facts as steps to further advancement. Computer formalisation is a big step forward for advancing certain parts of mathematics and a big philosophical advance which works together with, and not against, classical mathematics. The fact that human beings invented the printing press did not stop the development of writing, it enhanced it, and the same is true of any invention (weapons excluded) of our various civilisations. We should not fear that the computers are going to replace us- because they are not- we should welcome the new insights that they will make us have. 

We have already suggested elsewhere 
\cite{pluralism} that the solution of the crisis in the foundations might be to recognise the illusion of that ontological {\em the} in ``the foundations". The foundations are there to serve mathematics and not to bound its growth. If we need different foundations for different parts of mathematics, so be it, as long as we do not run into obvious contradictions. The foundations offered by set theory and category theory (which in our own view are different even if bi-interpretable) remain as important as ever, the univalent foundations have their place. 
\footnote{This view is easy to
take for a platonist, which the author happens to be, who anyway considers that all our mathematics is only a development towards the unreachable but existent world of real Mathematics and that the creation consists in being able to reach further towards that goal. How to do it is the essence of creativity. However, it does not seem to us that any of the known mathematical philosophy positions contradicts the evidence for pluralists fondations as represented by the above results. 

A `natural hero to logical pluralist', in the words of Greg Restall \cite{Restall} is Carnap, since,
still citing \cite{Restall},
`for him logic is of fundamental importance to philosophical inquiry, but nothing about this inquiry narrows the field of play to just one logic.'} 
Mathematics is also about intuition and it should be done by everybody in a way that maximises their own intuition, it may be ZFC for one and constructive mathematics for another. For example,
with our ZFC eyes,
we may pick from the article \cite{Bezem} the fact that there are $\Pi^0_1$-statements that cannot be proved in 
MLTT + Univalence and can be proved in ZFC. But, to a mathematician like Voevodsky it may be more striking that there are statements
about abstract category theory and abstract homotopy theory that are intuitively constructive
(and should not require strong proof theoretic means) but that become apparently
non effective, needing the Axiom of Choice, when formulated in ZFC.\footnote{We thank Thierry Coquand for the following example of
such a statement: a functor which is essentially surjective and full and faithful is an equivalence, as well as for the above comments.}.
When expressed in set theory, the proof requires the Axiom of Choice, which may be counter-intuitive. It would then be natural to look for a formalism
where such mathematical fact can be stated and proved in a more direct way, an intuition that it confirmed by the result about constructivity of univalence in \cite{Bezem},  which is the actual point of the paper from the point of view of its authors. 

In fact, a more careful reading of the history of logic than is within the scope of this article shows that various encodings of the foundations
have always coexisted.\footnote{Several readers of the first draft of the paper have commented that one should more didactically stress as a conclusion something that emerges as evidence from all stated above. So I will take the invitation to state what seems to me as an obvious conclusion, or more positively, a good beginning for future reflection. Pluralism is the only reasonable and honest choice in the foundations of mathematics today. Machines, including computers, bring forward human knowledge but do not replace it and it is by embracing what they have to offer that we are going to advance, not by positioning ourselves in camps in which ones fight with chalk and the others with computer chips.
 
Of course, we could explain and explain some more, but to quote Du{\v s}ko Pavlovi{\'c}, from his comments on this paper, which are so interesting and detailed that they could in themselves form another paper 
 \begin{quotation}
There are so many misunderstadings between good people for bad reasons, which it should be possible to clear. But the conundrum is: the more we try to explain, the more words we write, there is more chance for further misunderstandings.
I think the only way to resolve such a conundrum is honesty and trust. 
\end{quotation}}

Let us finish by a quote from the masters who wisely predict the fate of foundations of mathematics, in words that entirely represent our own view.
Bourbaki writes in the introduction to the first volume of \cite{Bourbakisettheory}:\footnote{The translation here is by the present author.}
\begin{quotation} We now know that, logically speaking, it is possible to derive almost all present-day mathematics from a unique source, set theory. By doing this we do not pretend to write a law in stone; maybe one day mathematicians will establish different reasoning which is not formalisable in the language that we adopt here and, according to some, recent progress in homology suggests that this day is not too far away. In this case one shall have to, if not, totally change the language, at least enlarge the syntax. It is the future of mathematics that will decide this.
 \end{quotation}

\begin{acknowledgement}
Many thanks to the organisers of the FOMUS conference in July 2016 for their invitation to give a talk and to participate in the panel discussion. I would also like to thank Peter Aczel, Andrej Bauer, Mark Bezem, Thierry Coquand, Laura Crossila, Deborah Kant, Angus Mcintyre, Marco Panza, 
Du{\v s}ko Pavlovi{\'c}, Michael Rathjen, Christian Rosendal
and Andr\'es Villaveces, as well as to the anonymous referee, for very useful and interesting discussions about various parts of this paper. My thanks equally go to the audiences in Paris, Nancy, Oxford, Teheran and Mexico City who have listened and contributed by their comments to the talk that accompanied the development of this paper.
\end{acknowledgement}

\bibliographystyle{plain}
\bibliography{../../bibliomaster}

\end{document}